\documentclass[11pt]{amsart}

\usepackage{amssymb,upgreek}

\usepackage{url}
\usepackage{bm}
\usepackage[left=1in,top=1in,right=1in,bottom=1in,head=.2in]{geometry}
\usepackage[titletoc]{appendix}

\usepackage[textsize=scriptsize]{todonotes}
\usepackage{color}
\usepackage{mathtools} 

\usepackage{fancyhdr}
\usepackage{mathrsfs}
\usepackage[cmtip,all]{xy}

\usepackage{hyperref}
\usepackage{enumerate}

\swapnumbers
\newtheorem{theorem}{Theorem}[section] 

\newtheorem*{theorem*}{Theorem}

\newtheorem*{fcthm*}{Finite Cork Theorem}
\newtheorem*{ccthm*}{Cork Consolidation Theorem}
\newtheorem*{thm1*}{Theorem 1}
\newtheorem*{thm2*}{Theorem 2}
\newtheorem*{mainthm*}{Main Theorem}
\newtheorem*{2lbthm*}{Generalized 4D Lightbulb Theorem (restated)}
\newtheorem*{icthms*}{Infinite Cork Theorems}
\newtheorem*{aclemma*}{\ac-Lemma}
\newtheorem*{mclemma*}{Multicork Lemma}
\newtheorem*{multicorktheorem*}{Multicork Theorem}
\newtheorem*{lemma*}{Lemma}
\newtheorem*{corollary*}{Corollary}

\newcommand{\thistheoremname}{}
\newtheorem{genericthm}[theorem]{\thistheoremname}

\theoremstyle{definition}

\newtheorem{remark}[theorem]{Remark}

\newtheorem*{remark*}{Remark}
\newtheorem*{definition*}{Definition}
\newtheorem*{remarks*}{Remarks}
\newtheorem*{addenda*}{Addenda}

\newcommand{\fig}[3]{\begin{figure}\includegraphics[height=#1pt]{#2}#3\end{figure}}

\newcommand{\bit}[1]{\textbf{\textit{#1}}} 

\newcommand{\sto}{\!\!\xymatrix@C=1em{{}\ar@{~>}[r]&{}}\!\!}

\newcommand{\del}{\partial}

\newcommand{\interior}{\textup{int}}
\newcommand{\ac}{\textup{AC}}

\newcommand{\items}{\begin{itemize}[leftmargin=25pt,rightmargin=5pt]
  \setlength\itemsep{2pt}}
\newcommand{\stopitems}{\end{itemize}}

\address{Princeton University,
Princeton, NJ 08544}
\email{hs25@princeton.edu} 

\begin{document}

\title{Duals of non-zero square}
\author{Hannah R. Schwartz}

\begin{abstract}
In this short note, for each non-zero integer $n$, we construct a $4$-manifold containing a smoothly concordant pair of spheres with a common dual of square $n$ but no automorphism carrying one sphere to the other. Our examples, besides showing that the square zero assumption on the dual is necessary in Gabai's and Schneiderman-Teichner's versions of the 4D Light Bulb Theorem, have the interesting feature that both the Freedman-Quinn and Kervaire-Milnor invariant of the pair of spheres vanishes. The proof gives a surprising application of results due to Akbulut-Matveyev and Auckly-Kim-Melvin-Ruberman pertaining to the well-known Mazur cork. 
\end{abstract}

\maketitle

\vskip-.4in
\vskip-.4in

\parskip 2pt

\setcounter{section}{-1}

\parskip 2pt

\section{Introduction and Motivation} \label{intro}

We work throughout in the smooth, oriented category. Begin by considering a pair of homotopic $2$-spheres $S$ and $T$ embedded in a smooth $4$-manifold $X$, with an embedded $2$-sphere $G \subset X$ intersecting both $S$ and $T$ transversally in a single point. Such a sphere is called a \bit{common dual} of $S$ and $T$. Recent work of Gabai \cite{dave:LBT} and Schneiderman-Teichner \cite{st} has completely characterized the conditions under which the spheres $S$ and $T$ are isotopic, so long as their common dual $G$ has square zero, i.e. a trivial normal bundle, in the $4$-manifold $X$. We call such a dual \bit{standard}, and \bit{non-standard} otherwise. The objective of this note is to show that the assumption of a standard dual is necessary one in \cite{dave:LBT} and \cite{st}. 

\begin{mainthm*} \label{diff2}
For each $n\not=0$, there exists a $4$-manifold $X_n$ containing smoothly concordant embedded spheres $S_n$ and $T_n$ with a common dual of square $n$ such that there is no automorphism of $X_n$ carrying one sphere to the other. 
\end{mainthm*} 

The proof of our Main Theorem gives a surprising application of well-studied $4$-dimensional objects called \bit{corks}: compact contractible 4-manifolds $C$ equipped with an orientation preserving diffeomorphism $h\colon\del C\to\del C$ . The study of corks was initially motivated by the fact that the \bit{cork twist} $X_{C,h} \ = \ (X-\interior(C))\,\cup_h \,C$ of an embedded cork $C \subset X$ is homeomorphic to $X$ by Freedman\,\cite{freedman:simply-connected}, but need not be diffeomorphic to $X$ by Akbulut~\cite{akbulut:contractible}. Such an embedding of a cork is called \bit{non-trivial}. Our construction builds upon examples given by Akbulut and Matveyev \cite{akbulut-matveyev:adjunction} of non-trivial embeddings of corks. 

\smallskip
\noindent
\bit{Acknowledgements.}
The author is grateful to both Dave Gabai and Peter Teichner for their encouragement to write up this result, and for their thoughtful comments and advice.

\section{Warm-up} \label{warmup}

The first example of a cork with a non-trivial embedding was produced by Akbulut in \cite{akbulut:contractible}. Now ubiquitous, the ``Akbulut-Mazur cork" $(W,\tau)$ consists of the Mazur manifold\footnote{Mazur's \cite{mazur} contractible $4$-manifolds are each built with a single $0$,$1$, and $2$-handle. They are not homeomorphic to the $4$-ball, but their products with the interval give the standard $5$-ball.} $W$ shown in Figure \ref{Mazur}, and $\tau$ the involution on its boundary induced by a rotation of $\pi$ around the indicated axis of symmetry. Many $4$-manifolds are now known to admit non-trivial embeddings of the Mazur cork; we outline one such embedding due to Akbulut and Matveyev \cite{akbulut-matveyev:adjunction} as a warm-up to the proof of the Main Theorem.

\fig{145}{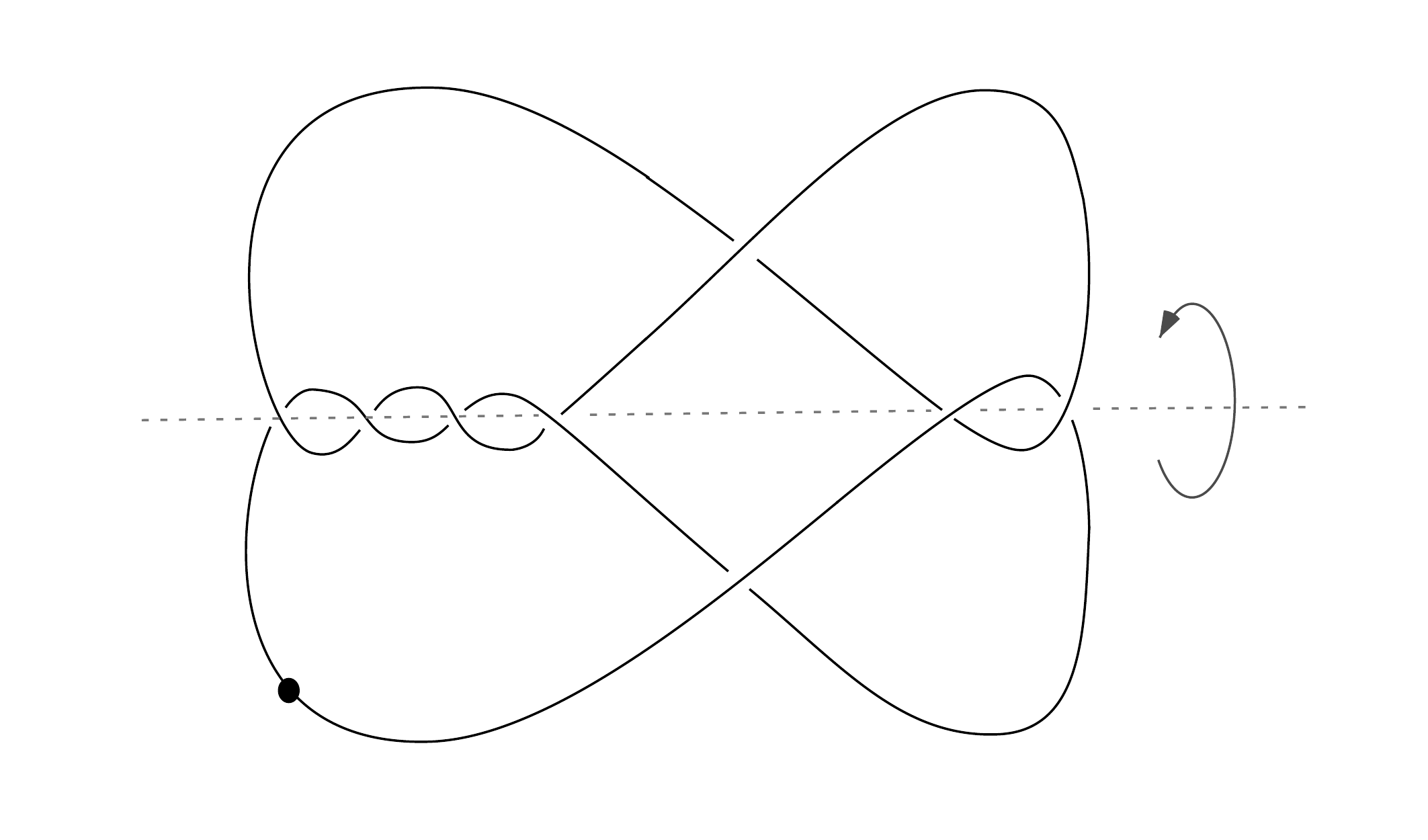}{
\put(-32,60){$\tau$}
\put(-214,100){$0$}
\caption{The Akbulut-Mazur cork $(W, \tau)$}
\label{Mazur}}

Let $X$ denote the compact 4-manifold shown on the left in Figure \ref{corktwist}, built from the Mazur manifold $W$ by adding a single 2-handle. Note that $X$ has a handlebody decomposition consisting of a single $1$-handle, and two $2$-handles each attached along knots in $S^1 \times S^2$ with framings less than their maximum Thurston-Bennequin numbers, as illustrated on the bottom right of Figure \ref{steinpic}. Therefore $X$ is a compact Stein domain\footnote{For a precise definition of what we mean by ``compact Stein domain", see \cite{akbulut-matveyev2}.} by a result of Eliashberg \cite{eliashberg:stein}; see also \cite{gompf} for more exposition.

\fig{145}{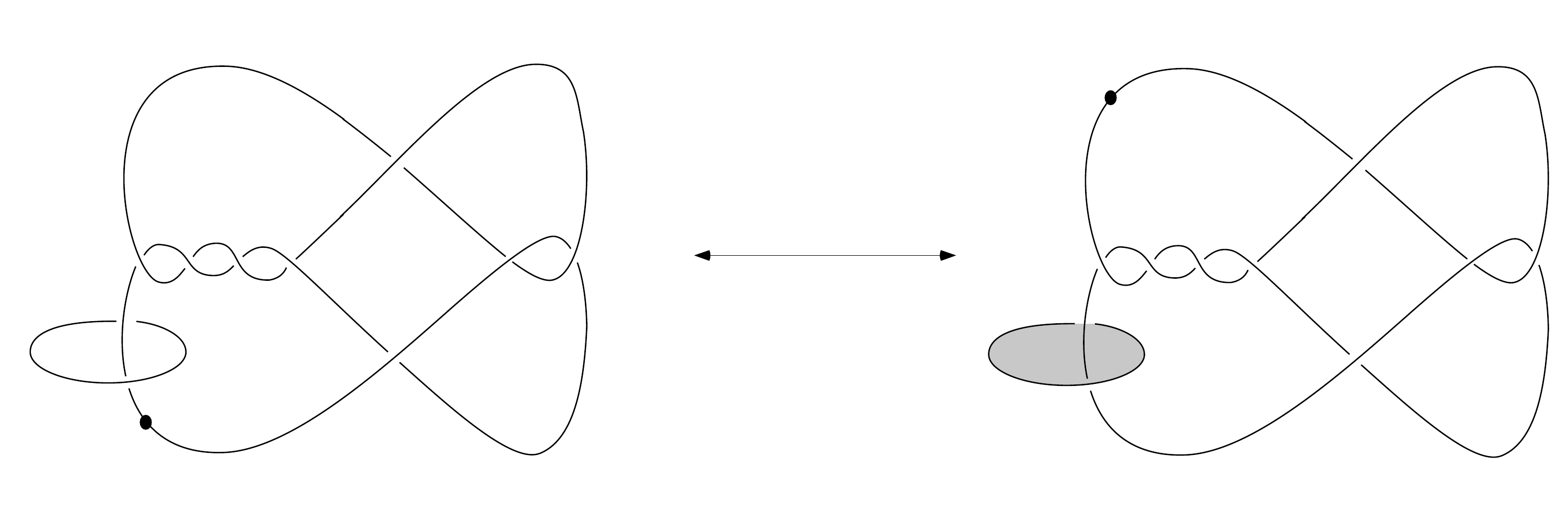}{
\put(0,43){$0$}
\put(-268,100){$0$}
\put(-177,42){$-1$}
\put(-149,43){\small{$D$}}
\put(-444,43){$-1$}
\put(-221,90){Mazur}
\put(-230,80){cork twist}
\caption{The manifold $X$ (left) and the cork twist $X_{W, \tau}$ (right)}
\label{corktwist}}

On the other hand, the cork twist $X_{W, \tau}$ contains an embedded 2-sphere of square $-1$, seen in the diagram for $X_{W, \tau}$ in Figure \ref{corktwist} as the union of the shaded disk $D$ and the core of the 2-handle attached along $\partial D$. Therefore $X_{W,\tau}$ must not be a compact Stein domain. This follows from a result due to Lisca and Mati\'c \cite{liscamatic} that compact Stein domains embed in minimal, closed K\"ahler surfaces, which contain no smoothly embedded 2-spheres of square $-1$. Therefore, $X$ and $X_{W,\tau}$ are not diffeomorphic. 

\fig{145}{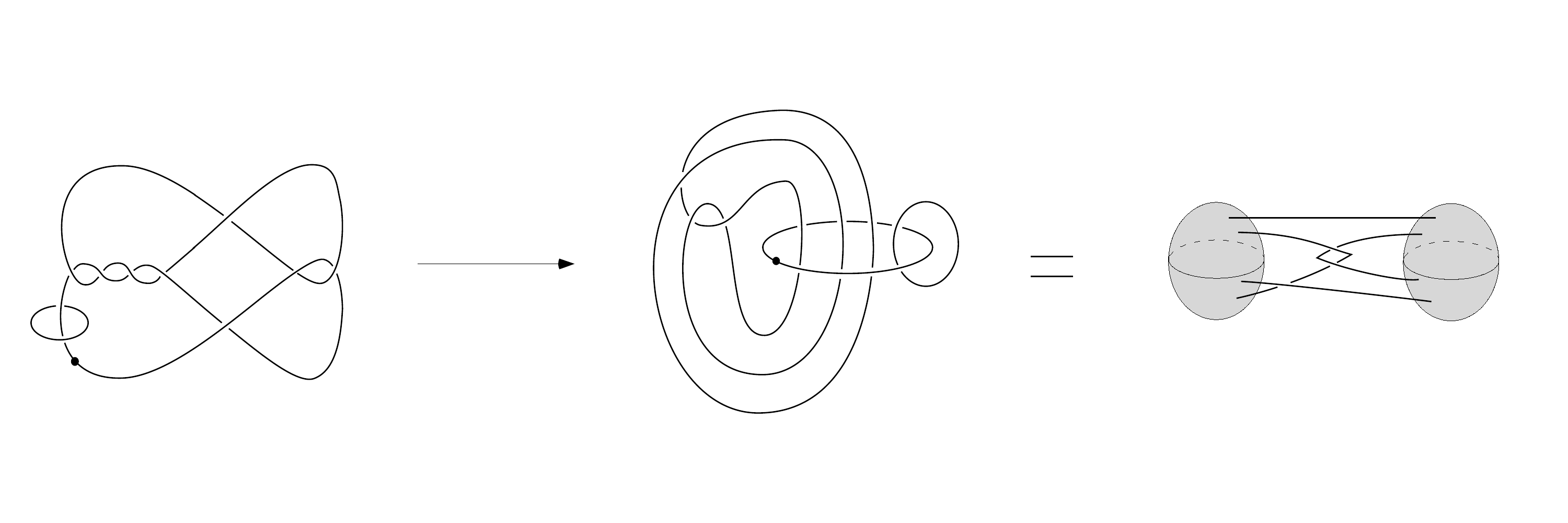}{
\put(-424,90){$0$}
\put(-260,48){$0$}
\put(-177,58){$-1$}
\put(-442,51){$-1$}
\put(-58,52){$0$}
\put(-63,88){$-1$}
\put(-316,76){Isotopy}
\caption{Identical handlebody structures for $X$, drawn with (left and middle) and without (right) the dotted circle notation for $1$-handles from \cite[Chapter I.2]{kirby:4-manifolds}. The Thurston-Bennequin framing of the attaching circle of each 2-handle is computed from the rightmost diagram using the usual formula (writhe) - (number of right cusps).}
\label{steinpic}}

\section{Main theorem} \label{main}

To contextualize our main result, we outline the previous results about common duals referred to in Section \ref{intro}. By Gabai \cite{dave:LBT} and Schneiderman-Teichner \cite{st}, the existence of a common standard dual for homotopic spheres $S, T \subset X$ guarantees a smooth isotopy between $S$ and $T$ whenever the \bit{Freedman-Quinn invariant}, a concordance invariant defined in \cite{fq}, of the pair $(S,T)$ vanishes. Recent work of Gabai \cite{dave:LBL} shows an analogous result holds for certain properly embedded disks with a common standard dual and vanishing \bit{Dax invariant}, an isotopy invariant of properly embedded disks recently formulated by Gabai in \cite{dave:LBL} using homotopy theoretic work of Dax \cite{dax} from the 70's. To guarantee even a smoothly embedded concordance between $S$ and $T$ when their common dual is non-standard, it is also required that their \bit{Kervaire-Milnor invariant}, defined by Stong in \cite{stong}, vanishes. 

\begin{remark} Recently, Klug-Miller \cite[Example 7.2]{mm} showed that it is necessary that the dual have square zero for Gabai \cite{dave:LBT} and Schneiderman-Teichner \cite{st} to achieve an isotopy, by presenting a pair of spheres whose common dual of square $+1$, with vanishing Freedman-Quinn invariant but non-vanishing Kervaire-Milnor invariant. On the other hand, for each $n \not =0$, the Main Theorem gives examples of pairs of spheres with dual of square $n$ whose Freedman-Quinn invariant and Kervaire-Milnor invariants vanish, but that are not related by any automorphism of the ambient $4$-manifold. Such an automorphism always exists for spheres when the common dual is standard by \cite[Lemma 2.3]{me}, since in this case the common dual can be surgered (Gabai remarks after \cite[Theorem 0.8]{dave:LBL} that by a similar proof, this also holds for properly embedded disks with a common standard dual). 
\end{remark}

\smallskip
\noindent \emph{Proof of Main Theorem.} For $n\leq -1$, consider the $4$-manifold $X_n$ pictured in Figure \ref{nonzero}. Since $X_n$ is simply-connected, the spheres $S_n$ and $T_n$ are not only homologous, but also homotopic. It is also immediate that the both the Freedman-Quinn and Kervaire-Milnor invariants of the pair $(S_n, T_n)$ vanish, since these invariants are elements of $H_1(X_n; \mathbb{Z}_2)$ and a quotient of $\mathbb{Z}[\pi_1(X)]$ respectively, which are both trivial in this case. Let $R_n$ denote the sphere of square $n$ gotten by capping off the red disk in Figure \ref{nonzero} with the core of the $2$-handle attached with framing $n$ along its boundary. The sphere $R_n$ is dual to both $S_n$ and $T_n$, since $S_n$ and $T_n$ each pass once (geometrically) over the $2$-handle with framing $1$ in the topmost diagram of Figure \ref{nonzero}. Therefore, by \cite{fq} and \cite{stong}, the spheres $S_n$ and $T_n$ are smoothly concordant in $X_n \times I$. 

\fig{340}{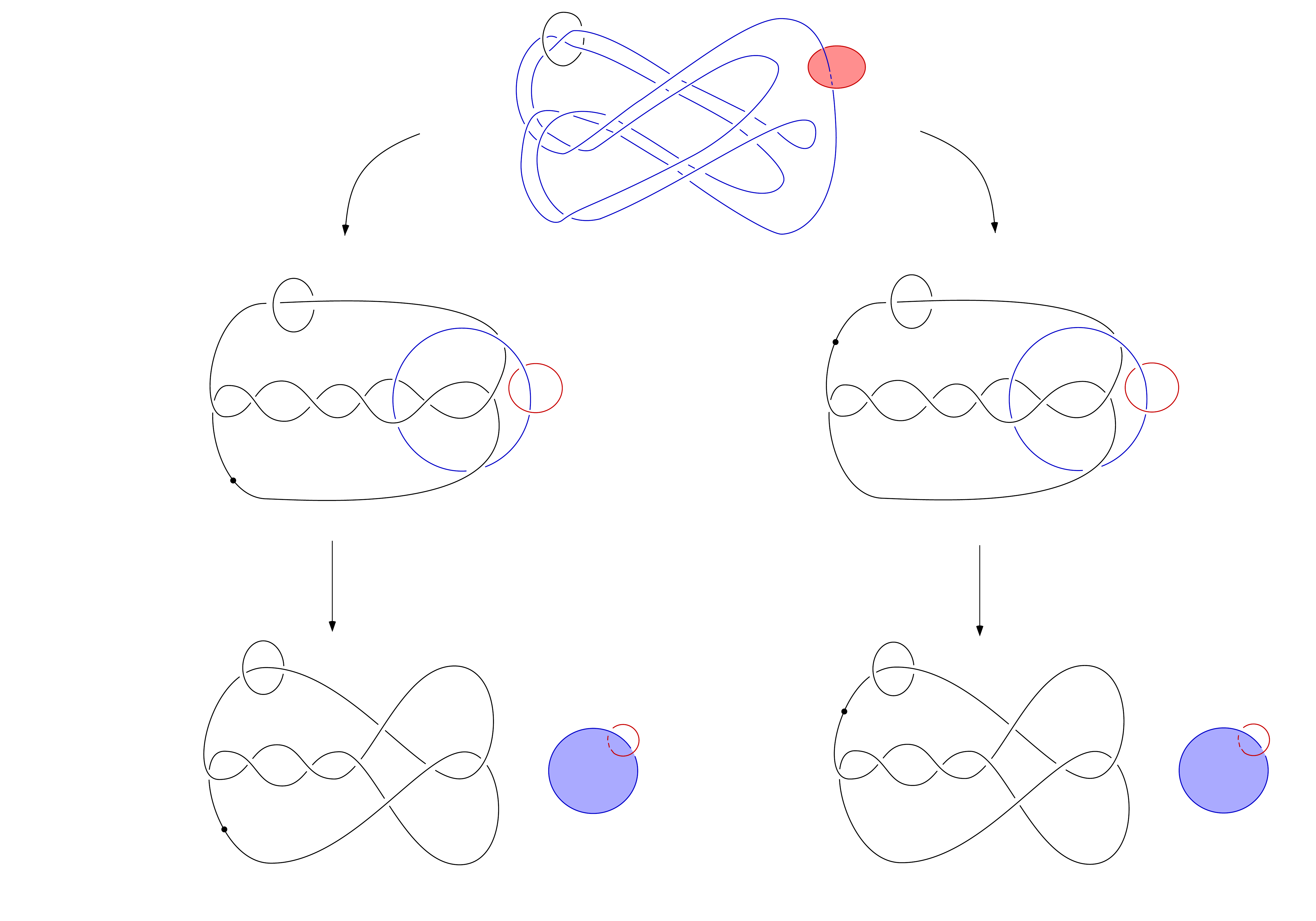}{
\put(-182,32){$0$}
\put(-418,75){$0$}
\put(-188,167){$0$}
\put(-417,209){$0$}
\put(-346,209){$1$}
\put(-118,210){$1$}
\put(-305,310){$1$}
\put(-285,338){$-1$}
\put(-382,240){$-1$}
\put(-154,241){$-1$}
\put(-160,106){$-1$}
\put(-390,105){$-1$}
\put(-283,30){$1$}
\put(-45,30){$1$}
\put(-248,64){$n$}
\put(-15,64){$n$}
\put(-47,193){$n$}
\put(-275,192){$n$}
\put(-163,312){$n$}
\put(-271,50){$S_n$}
\put(-39,50){$T_n$}
\put(-113,120){Handleslide}
\put(-428,120){Handleslide}
\put(-100,290){Introduce}
\put(-110,270){cancelling 1/2 pair}
\put(-430,290){Introduce}
\put(-460,270){cancelling 1/2 pair}
\caption{The spheres $S_n$ and $T_n$ (blue) in $X_n$, with their common dual (red).}
\label{nonzero}}

The manifold $X_n$ contains Akbulut and Matveyev's manifold $X$ \cite{akbulut-matveyev:adjunction} discussed in Section \ref{warmup}. To show that there is no automorphism of $X_n$ carrying $S_n$ to $T_n$, we use an argument similar to one of Auckly-Kim-Melvin-Ruberman \cite[Theorem A]{akmr:isotopy}; see in particular Figure $18$ of their paper. For, blowing down $S_n$ gives the bottom left manifold of Figure \ref{nonzero211}, which is not Stein since it contains an embedded sphere of square $-1$, as in the argument from Section \ref{warmup}. On the other hand, blowing down $T_n$ gives the bottom right manifold of Figure \ref{nonzero211}, which is Stein whenever $n \leq -1$ by \cite{eliashberg:stein}, since all $2$-handles are attached along Legendrian knots whose framings are strictly less than their Thurston-Bennequin numbers. 

As the manifolds that result from blowing down $S_n$ and $T_n$ are not diffeomorphic, there can be no automorphism of $X_n$ carrying one sphere to the other when $n \leq -1$. The result therefore also holds for $n \geq 1$, setting $X_n=-X_{-n}$ and considering the spheres $S_n, T_n \subset X_n$ that are the images of the spheres $S_{-n}, T_{-n} \subset X_{-n}$ under the (orientation reversing) identity map. \qed 

\fig{310}{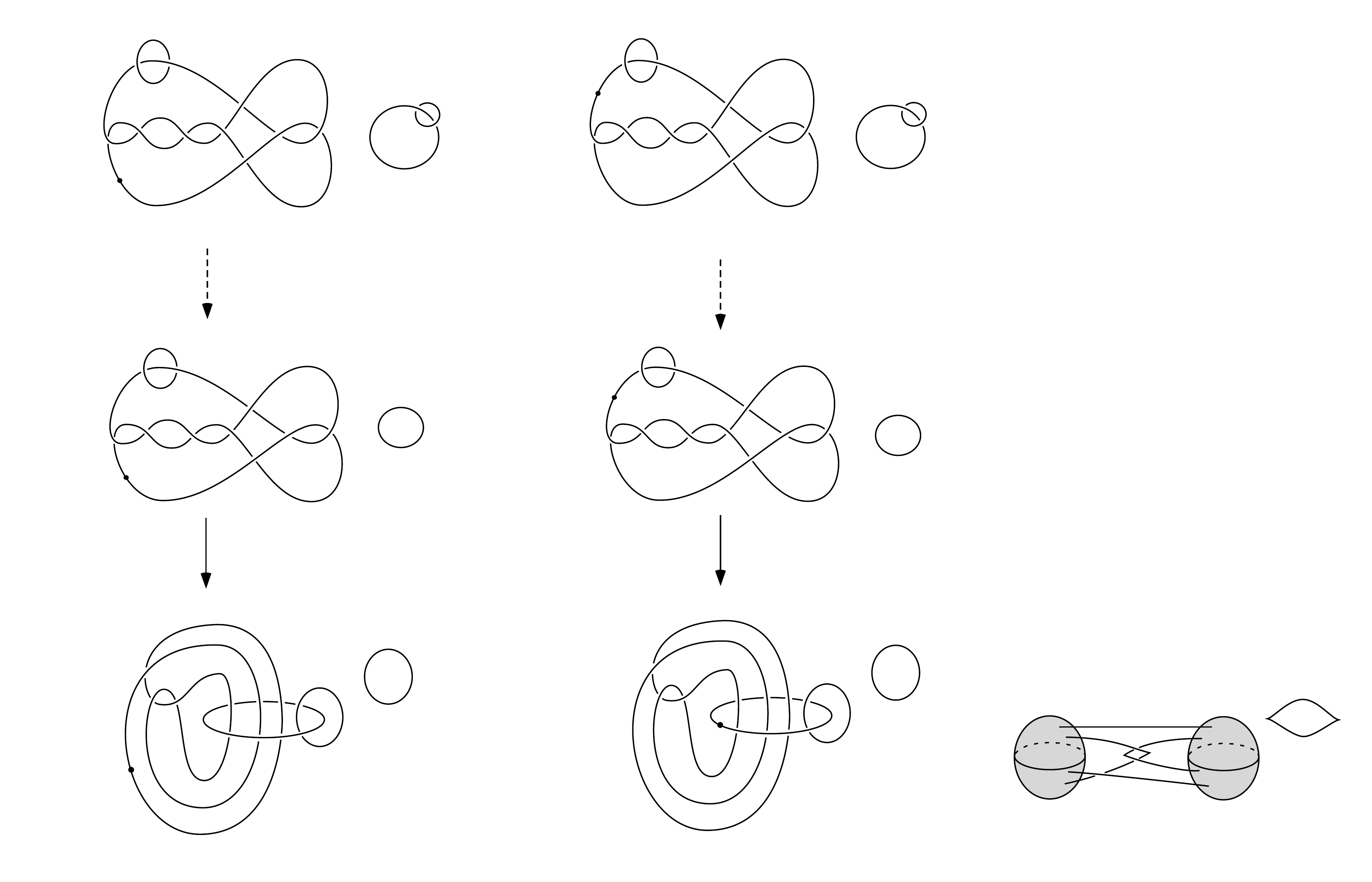}{
\put(-12,63){\small{$n-1$}}
\put(-155,40){$=$}
\put(-90,58){\small{$-1$}}
\put(-80,25){\small{$0$}}
\put(-155,270){\small{$n$}}
\put(-325,270){\small{$n$}}
\put(-438,300){\small{$-1$}}
\put(-266,300){\small{$-1$}}
\put(-437,192){\small{$-1$}}
\put(-260,192){\small{$-1$}}
\put(-378,37){\small{$-1$}}
\put(-200,38){\small{$-1$}}
\put(-452,273){\small{$0$}}
\put(-280,245){\small{$0$}}
\put(-450,167){\small{$0$}}
\put(-275,140){\small{$0$}}
\put(-413,44){\small{$0$}}
\put(-211,22){\small{$0$}}
\put(-172,240){\small{$1$}}
\put(-343,240){\small{$1$}}
\put(-180,86){\small{$n-1$}}
\put(-359,85){\small{$n-1$}}
\put(-355,170){\small{$n-1$}}
\put(-180,168){\small{$n-1$}}
\put(-225,112){\small{Isotopy}}
\put(-404,112){\small{Isotopy}}
\put(-225,207){\small{Blow down $T_n$}}
\put(-404,209){\small{Blow down $S_n$}}
\caption{Blowing down the spheres $S_n$ and $T_n$}
\label{nonzero211}}

\end{document}